\documentclass[10pt]{amsart}
\usepackage{amssymb,latexsym}

\newtheorem{theorem}{Theorem}

\newtheorem*{TA}{Theorem A}
\newtheorem*{TB}{Theorem B}
\newtheorem*{TC}{Theorem C}
\newtheorem*{TD}{Theorem D}

\newtheorem{lemma}{Lemma}

\theoremstyle{definition}

\theoremstyle{remark}

\numberwithin{equation}{section}

\newcommand{\sgn}{\text{sgn}}

\newcommand{\R}{\mathbb R}

\begin{document}
\title[Persistence properties]{
On persistence properties in fractional weighted spaces}
\author{G. Fonseca}
\address[G. Fonseca]{Departamento  de Matem\'aticas\\
Universidad Nacional de Colombia\\ Bogot\'a\\Colombia}
\email{gefonsecab@unal.edu.co}
\author{F.  Linares}
\address[F. Linares]{IMPA\\
Instituto Matem\'atica Pura e Aplicada\\
Estrada Dona Castorina 110\\
22460-320, Rio de Janeiro, RJ\\Brazil}
\email{linares@impa.br}
\author{G. Ponce}
\address[G. Ponce]{Department of Mathematics\\ South Hall, Room 6607\\ University of California,
Santa Barbara\\CA 93106, USA}
\email{ponce@math.ucsb.edu}
\thanks{The third author is supported  by  NFS grant DMS-1101499}
\begin{abstract} In this work we derive a point-wise formula that will allows us to study the well-posedness of initial value problem associated to nonlinear dispersive equations
in fractional weighted Sobolev spaces $H^s(\R)\cap L^2(|x|^{2r}dx)$, $s, r \in \R$. As an application of this formula we will study local and global well posedness of  
the $k$-generalized Korteweg-de Vries equation  in these weighted Sobolev spaces.
\end{abstract}
\maketitle

\section{Introduction}\label{S:1}
In this work we are concerned with persistence properties of solutions of the initial value problem (IVP) associated  to nonlinear dispersive equations
in fractional weighted spaces. More precisely, if we define the weighted Sobolev spaces
\begin{equation}
\label{spaceZ}
Z_{s,r} = H^s(\R) \cap L^2(|x|^{2r}dx),\,\,\, s, r\in \R ,
\end{equation}
we would like to prove that for data in the function space the associated IVP is locally or globally well-posed. We will follow the notion of  well posedness given  in \cite{Ka1}:  
 the  IVP is said to be locally well posed (LWP) in the function space $X$ if for each $u_0\in X$ there exist $T>0$ and a unique solution $u\in C([-T,T]:X)\cap\dots=Y_T$ of 
the equation, with the map data $\to$ solution being locally continuous from $X$ to $ Y_T$. 

This notion of LWP includes the \lq\lq  persistence" property, i.e. the solution describes a continuous curve on $X$. In particular, this implies that the solution flow of the considered 
equation defines a dynamical system in $X$.   If $T $ can be taken arbitrarily large, then the IVP  is said to be globally well posed (GWP).

To present our main result and give some applications  we will use as example the IVP associated to the $k$-generalized Korteweg-de Vries
equation, 
\begin{equation}\label{kgKdV}
\begin{cases}
\partial_tu +\partial_x^3 u+ u^k\partial_xu=0,\;\;\;\;t,\,x\in\R,\;\;k\in\mathbb Z^+,\\
u(x,0)=u_0(x).
\end{cases}
\end{equation}
However, the main result is quite general as we will comment below.

Concerning LWP in the weighted spaces $Z_{s,r}$ defined in \eqref{spaceZ} T. Kato  \cite{Ka1}  showed that 
 persistent properties hold for solutions of the IVP \eqref{kgKdV} for any $m\in\mathbb Z^+$ in
$$
Z_{s,m}=H^s(\R)\cap L^2(\,|x|^{2m}),\;\;\;\;\;s\geq 2m,\;\;\;\;\;m=1,2,\dots
$$
More precisely:

\begin{TA} $($\cite{Ka1}$)$
Let $m\in\mathbb Z^+$.  Let $u\in C([-T,T]:H^{s}(\R))\cap \dots$ with $s\geq 2m$ be the solution of the IVP  \eqref{kgKdV}. If  $u(x,0)=u_0(x)\in L^2(|x|^{2m}dx)$, then
$$
u\in C([-T,T]:Z_{s,m}).
$$
\end{TA}

The proof of Theorem A in  \cite{Ka1}  relies on the commutative property of the operators
\begin{equation}
\label{op-ga}
\Gamma =x-3t\partial_x^2,\;\;\;\;\;\;\;\mathcal L=\partial_t+\partial_x^3,\;\;\;\;\;\text{so}\;\;\;\;[\Gamma;\mathcal L]=0.
\end{equation}
 In particular, if $\,\{U(t) \,:\,t\in\R\}$ denotes the unitary group of operators describing the solution of the linear IVP 
 \begin{equation}
\label{lKdV}
\partial_tv +\partial_x^3 v=0,\;\;\;\;t,\,x\in\R,\;\;\;\;\;\;\;\; v(x,0)=v_0(x),
\end{equation}
i.e. 
\begin{equation}\label{de-gr}
U(t)v_0(x)= (e^{it\xi^3} \widehat v_0)^{\lor}(x),
\end{equation}
 then from  \eqref{op-ga} one has that
 \begin{equation}
 \label{formula1}
 \begin{aligned}
& x \,U(t) v_0(x) = U(t)(xv_0)(x) +3t U(t) (\partial_x^2 v_0)(x),\\
&\;\;\;\;\text{i.e.}\;\;\;\;\\
&\Gamma U(t)v_0(x)=U(t)(xv_0)(x).
\end{aligned}
 \end{equation}
 
The form of the operator $\Gamma$ suggests  that one should expect persistence in $Z_{s,r}$ only if $s\geq 2r$.  Thus in order to treat
fractional powers of $x$ (or $|x|$) we would like to have an identity in the same spirit as \eqref{formula1}. This is what our main result guarantees.
More precisely we shall prove:
 \vskip.1in
 \begin{theorem} \label{point-wise}
Let $\alpha\in (0,1)$ and $\{U(t)\,:\,t\in\R\}$ be the unitary group of operators defined in \eqref{de-gr}.   If
\begin{equation} \label{hyp1}
u_0\in  Z_{2\alpha,\alpha}=H^{2\alpha}(\R)\cap L^2(|x|^{2\alpha}dx),
\end{equation}
then for all $t\in \R$ and for almost every $x\in \R$
\begin{equation}
\label{d20}
|x|^{\alpha} U(t) u_0(x) = U(t) (|x|^{\alpha} u_0)(x) + U(t)\{\Phi_{t,\alpha}(\widehat u_0)(\xi)\}^{\lor}(x)
\end{equation}
 with
 \begin{equation}
\label{d20-norm}
\| \{ \Phi_{t,\alpha}(\widehat u_0)(\xi)\}^{\lor}\|_2\leq c(1+|t|)(\|u_0\|_2 + \| D^{2\alpha}u_0\|_2).
\end{equation}
Moreover, if in addition to \eqref{hyp1} one has  that for $\beta\in (0,\alpha)$
\begin{equation}
\label{hyp2}
D^{\beta} (|x|^{\alpha}u_0)\in  L^{2}(\R)\;\;\;\;\;\text{and}\;\;\;\;\;u_0\in H^{\beta+2\alpha}(\R),
\end{equation}
 then for all $t\in \R$ and for almost every $x\in \R$
\begin{equation}
\label{d21}
\begin{aligned}
&D^{\beta} (|x|^{\alpha} U(t) u_0)(x) \\
\\
&= U(t) (D^{\beta}|x|^{\alpha} u_0)(x) + U(t)(D^{\beta} (\{ \Phi_{t,\alpha}(\widehat u_0)(\xi)\}^{\lor}))(x)
\end{aligned}
\end{equation}
  with
 \begin{equation}
\label{d21-norm}
\|D^{\beta}(\{ \Phi_{t,\alpha}(\widehat u_0)(\xi)\}^{\lor})\|_2\leq c(1+|t|)(\|u_0\|_2 + \| D^{\beta+2 \alpha}u_0\|_2).
\end{equation}
 \end{theorem}

\underline{Remark:}
The identities \eqref{d20}-\eqref{d20-norm} can be seen as an extension of \eqref{formula1} for fractional weights.
 As it will be remarked below the result in Theorem \ref{point-wise} can be adapted to general groups describing the solution 
 of the linear part of a dispersive equation.

The  proof of Theorem \ref{point-wise}  will be based on a characterization of the generalized Sobolev space
 \begin{equation}
 \label{def0}
 L^{\alpha,p}(\R^n)= (1-\Delta)^{-\alpha/2} L^{p}(\R^n),\;\;\;\;\alpha\in (0,2),\;\;p\in(1,\infty),
 \end{equation}
 due to E. M. Stein \cite{St1} (see Theorem D below).

As we mentioned above as an application of our main result we will study persistence properties of solutions of the initial value problems (IVP) associated to  the $k$-generalized 
Korteweg-de Vries ($k$-gKdV) equation \eqref{kgKdV} in weighted Sobolev spaces
 \begin{equation}
 \label{w-s}
 Z_{s,r}\equiv H^s(\R)\cap L^2(\,|x|^{2r}),\;\;\;\;\;s\in \R,\;\,r\geq 0.
 \end{equation}

 We shall be mainly concerned with  the modified Korteweg-de Vries  (mKdV) equation, i.e. $k=2$ in \eqref{kgKdV}.
 In \cite{KPV1} Kenig, Ponce and  Vega showed that the IVP  \eqref{kgKdV} with $k=2$ is locally well posed in
 $$
 \dot H^{1/4}(\R)=(-\partial_x^2)^{-1/8} L^2(\R)\supset H^{1/4}(\R)=J^{-1/4}L^2(\R)=(1-\partial_x^2)^{-1/8}L^2(\R).
 $$
 More precisely, the following result was established in \cite{KPV1}:

\begin{TB}$($\cite{KPV1}$)$
For any $u_0\in\dot H^{1/4}(\R)$ there exist
\begin{equation}
\label{time-es}
T=T(\|D_x^{1/4}u_0\|_2)\sim \|D_x^{1/4}u_0\|_2^{-4},
\end{equation}
and a unique solution $u(t)$ of the IVP \eqref{kgKdV} with $k=2$ such that
\begin{equation}
\label{class1}
\begin{aligned}
& \;\;\;\;u\in C([-T,T]:\dot H^{1/4}(\R)),\\
&\text{and}\\
&\;\;\;\;\|D_x^{1/4} \,\partial_x u\|_{L_x^\infty L_T^2}+\|\partial_x u\|_{L_x^{20}L_T^{5/2}}+\|D_x^{1/4}u\|_{L_x^5 L_T^{10}}+ \|u\|_{L_x^4 L_T^\infty}<\infty.
\end{aligned}
\end{equation}

For any $T'\in(0,T)$ there exists a neighborhood $V$ of $u_0$ in
$\dot H^{1/4}(\R)$ such that the map data $\to$ solution $\widetilde u_0\to\widetilde
u(t)$ from $V$ into the class defined by \eqref{class1} with $T'$
instead of $T$ is smooth.

Moreover, if in addition $u_0\in H^s(\R)$ with $\,s\geq 1/4$, then the solution 
$$
u \in C([-T,T]: H^s(\R)),
$$
and 
$$
\|D_x^{s} \,\partial_x u\|_{L_x^\infty L_T^2}+ \|J^{s-1/4}_x\partial_xu\|_{L_x^{20}L^{5/2}_T}+\|J_x^{s}
u \|_{L_x^5 L_T^{10}}<\infty.
$$
\end{TB}

\underline{Remarks}: (a) The fact that the map data $\to$ solution is smooth is a direct consequence of the proof
of  Theorem B, based on the contraction principle, and the implicit function theorem. The estimate for the length 
of the time interval of existence \eqref{time-es} is  inside the proof in \cite{KPV1} (which is partially reproduced in the proof of Theorem \ref{th4} below)
or can also be obtained by a scaling argument. 

(b) It was shown in  \cite{KPV3} and \cite{CCT} that in an appropriate sense the value $1/4$ in Theorem B  is optimal.

(c) In \cite{CKSTT} Colliander, Keel, Staffilani, Takaoka, and Tao showed that this LWP extends to 
a GWP if $s>1/4$. The GWP for the limiting case $s=1/4$ was established by Guo \cite{Gu} and
Kishimoto \cite{Ki1}. 

(d)  We recall the best known  LWP and GWP results in $H^s(\R)$  for the IVP \eqref{kgKdV} with $k\neq 2$:

 - for  $k=1$  LWP is known for  $s\geq -3/4$ (see  \cite{KPV2} for the case $s>-3/4$  and \cite{CCT}, \cite{Gu} and \cite{Ki1} for the limiting case $s=-3/4$), and  GWP is known for $s\geq -3/4$ (see  \cite{CKSTT} for the case $s>-3/4$ and  \cite{Gu} and \cite{Ki1} for the limiting case $s=-3/4$),
 
 - for  $k=3$ LWP is known for $s\geq -1/6$ (see \cite{AG} for the case $s>-1/6$ and \cite{Tao2} for the limiting case $s=-1/6$)
 and GWP is known for $s> - 1/42$ (see \cite{AGP}),
 
 -for $k\geq 4$ LWP is known for $s\geq (k-4)/2k$ (see \cite{KPV1}). In \cite{MM} for the case $k=4$ it is shown 
 that there exist local smooth solutions which develop singularities in finite time.

\begin{TC}$($\cite{Gu}, \cite{Ki1}$)$ 
 Let $u_0\in H^s(\R)$ with $s\geq 1/4$. Then for  any  $T^*>0$  the IVP  \eqref{kgKdV}
with $k=2$ has a unique solution
\begin{equation}
\label{class7}
u\in C([-T^*, T^*]: H^s(\R))\cap\dots
\end{equation}
\end{TC}

 Remark: (a) The proof of Theorem C relies on the so called \lq\lq I-method" introduced in \cite{CKSTT0}, on the Miura transformation \cite{Mi}, and on sharp  LWP for the Korteweg-de Vries (KdV) $k=1$ in \eqref{kgKdV}. This optimal LWP result for the KdV requires the use of the so called Bourgain spaces $X_{s,b}$, introduced in the context of non-linear dispersive equations in  \cite{Bo1}. Consequently, the precise description of the class in \eqref{class7} involves those spaces.

(b) In \cite{JN2} for the case of the mKdV, J. Nahas extended  locally the result in Theorem C to the optimal range of the parameter $s, r$ accordingly to Theorem A and \eqref{op-ga}, i.e. $s\geq 1/4$ and
$s\geq 2r>0$. Also in  \cite{JN2} for the case $k\geq 4 $ in \eqref{kgKdV} 
Theorem C was extended to the optimal range $s\geq (k-4)/4k$ and $s\geq 2r>0$.

Our second result  gives a significantly simplified proof and slightly stronger version of these results.
We shall concentrate in the case of the mKdV equation $k=2$ in \eqref{kgKdV}.

\begin{theorem}
\label{th4}

Let  $u\in C([-T,T]:\dot H^{1/4}(\R))$  denote the solution of the IVP  \eqref{kgKdV} with $k=2$ provided by Theorem A. If $\,u_0,\,|x|^{r}u_0\in L^2(\R)$ with $r\in(0,1/8]$, then 
\begin{equation}
\label{class2}
\;\;u\in C([-T,T]: Z_{1/4,r}).
\end{equation}

For any $T'\in(0,T)$ there exists a neighborhood $V$ of $u_0$ in
$ H^{1/4}(\R)\cap L^2(|x|^{2r}dx)$ such that the map $\widetilde u_0\to\widetilde
u(t)$ from $V$ into the class defined by \eqref{class1} and \eqref{class2} with $T'$
instead of $T$ is smooth.

Moreover, if in addition $u_0\in Z_{s,r}$ with $s>1/4$ and  $s\geq 2r$, then the solution 
\begin{equation}
\label{1.11}
u \in C([-T,T]: Z_{s,r}).
\end{equation}
\end{theorem}

\underline{Remarks}:  (a) We observe that Theorem \ref{th4} guarantees that the  persistent property in the weighted space 
$Z_{s,r}$ holds in the same time interval $[-T,T]$  given by Theorem A, where
$T$ depends only on $\,\|D^{1/4}_xu_0\|_2$ (see \eqref{time-est}). 

(b) It was established  in \cite{ILP} that the condition $s\geq 2r$ in Theorem \ref{th4} is optimal.
More preciesely, \eqref{1.11} can hold only if $s\geq 2r$.

(c) Roughly, in \cite{GT} Ginibre and Tsutsumi obtained 
results concerning the uniqueness and existence 
(in an appropriate class) of local solutions of the IVP \eqref{kgKdV} with $k=2$ and data $u_0$ in the weighted space
$L^2((1+|x|)^{1/4} dx)$. Theorem \ref{th4} shows that for data
 $u_0\in Z_{1/4, 1/8}$ the solution provided by Theorem A  and that obtained in \cite{GT} agree.

(d) As in \cite{JN2} the result in Theorem \ref{th4} extends to the local solutions of the IVP \eqref{kgKdV} with $k\geq 4$ in the optimal range of the parameters $s, r$ accordingly to remark (a) after Theorem  C, i.e.  $s\geq 2r>0$ with $s\geq (k-4)/2k$. This will be clear from our proof of Theorem \ref{th4} given below, so we omit the details.  For the cases $k=1$ and $k=3$ a weaker version of these  results was proven in \cite{NP}. The main difference between the cases
$k=2,4,5,...$ and $k=1,3$ is that for the latter the \lq\lq optimal"  well-posedness results are based on the spaces $X_{s,b}$ which make  fractional weights difficult to handle. 

As a consequence of Theorem B and our proof of Theorem \ref{th4} we obtain the following global version of Theorem \ref{th4}:

\begin{theorem}
\label{th5}

Let $s\geq 1/4$ and $\,T^*>0$. If $u_0\in Z_{s,r}$ with $s\geq2r>0$, then the solution 
$u$ of the IVP \eqref{kgKdV} with $k=2$ provided by Theorem \ref{th4} extends to the time interval 
$[- T^*, T^*] $ with
$$
u\in C([-T^*, T^*]\,:\,Z_{s,r}).
$$
\end{theorem}

The paper is organized as follows. The proof of Theorem \ref{point-wise} will be given in Section 2.  In Section 3 we will present the proofs of Theorem \ref{th4} and Theorem \ref{th5}.

\section{Proof of Theorem \ref{point-wise} }

Next we turn our attention to the proof of Theorem \ref{point-wise}. 
We shall start with a characterization of the Sobolev space
 \begin{equation}
 \label{def00}
 L^{\alpha,p}(\R^n)= (1-\Delta)^{-\alpha/2} L^{p}(\R^n),\;\;\;\;\alpha\in (0,2),\;\;p\in(1,\infty),
 \end{equation}
due to E. M. Stein \cite{St1}. For $\alpha \in (0,2)$ define

 \begin{equation}
 \label{def1}
 D_{\alpha}f(x) = \lim_{\epsilon \to 0} \frac{1}{c_{\alpha}} \,\int_{|y| \geq \epsilon}  \frac{f(x+y)-f(x)}{|y|^{n+\alpha}}dy,
 \end{equation}
where $\;c_{\alpha}=\pi^{n/2}\,2^{-\alpha}\,\Gamma (-\alpha/2)/\Gamma((n+2)/2)$.

  As it was remarked in  \cite{St1} for appropriate $f$, for example $f\in \mathcal S(\R^n)$, one has
 \begin{equation}
 \label{pro1}
\widehat{ D_{\alpha}f}(\xi) = \widehat{D^{\alpha} f}(\xi)\equiv  |\xi|^{\alpha} \,\widehat f(\xi).
\end{equation}

 The following  result concerning  the $L^{\alpha,p}(\R^n)= (1-\Delta)^{\alpha/2} L^{p}(\R^n)$ spaces
 was established  in \cite{St1},
\begin{TD} \label{Theorem4}
Let $\alpha\in (0,2)$ and $p\in(1,\infty)$. Then $f\in  L^{\alpha,p}(\R^n)$ if and only if
\begin{equation}
\label{d1}
\begin{cases}
&\;(a)\;\, f\in L^p(\R^n),\\
\\
&\;(b)\;\;D_{\alpha} f\in L^p(\R^n),\;\;\;\;\;\;\,\;\;\;\;\;\;\;\;\;\;\;\;(D_{\alpha} f (x)\;\;\text{defined in \eqref{def1}}),
\end{cases}
\end{equation}
with
 \begin{equation}\label{d1-norm}
\|f\|_{\alpha,p}= \|(1-\Delta)^{\alpha/2} f\|_p\simeq \|f\|_p+\|D_{\alpha}    f\|_p\simeq \|f\|_p+\|\,D^{\alpha}       f\|_p.
\end{equation}
 \end{TD}

 Notice that if  $ f,\, fg\in L^{\alpha,p}(\R^n)$ and $ g\in L^{\infty}(\R^n)\cap C^2(\R^n)$ one has

 \begin{equation}
\label{d2}
\begin{aligned}
D_{\alpha}(fg)(x) & = \lim_{\epsilon \to 0} \frac{1}{c_{\alpha}} \,\int_{|y| \geq \epsilon}  \frac{f(x+y)\,g(x+y)-f(x)\,g(x)}{|y|^{n+\alpha}}dy\\
\\
&= \lim_{\epsilon \to 0} \frac{1}{c_{\alpha}} \int_{|y| \geq \epsilon}  g(x)\frac{f(x+y)-f(x)}{|y|^{n+\alpha}}dy\\
\\
&\;\;\;+ \lim_{\epsilon \to 0} \frac{1}{c_{\alpha}} \int_{|y| \geq \epsilon}  \frac{(g(x+y)-g(x))f(x+y)}{|y|^{n+\alpha}}dy\\
\\
&= g(x)\,D_{\alpha}f(x) + \Lambda_{\alpha}((g(\cdot+y)-g(\cdot))f(\cdot+y))(x).
\end{aligned}
\end{equation}

  In particular, if $ g(x) =e^{i\phi(x)}$, then
  \begin{equation}
\label{d3}
\begin{aligned}
  &\Lambda_{\alpha}((g(\cdot+y)-g(\cdot))f(\cdot+y))(x)\\
  \\
  & =
   \lim_{\epsilon \to 0} \frac{1}{c_{\alpha}} \int_{|y| \geq \epsilon}  \frac{(g(x+y)-g(x))f(x+y)}{|y|^{n+\alpha}}dy\\
   \\
  &= e^{i\phi(x)}\, \lim_{\epsilon \to 0}\, \frac{1}{c_{\alpha}}\int_{|y| \geq \epsilon} \frac{e^{i(\phi(x+y)-\phi(x))}-1}{|y|^{n+\alpha}}\,f(x+y) dy.
  \end{aligned}
  \end{equation}

  Thus,  one gets the identity
  \begin{equation}
\label{d4}
  D_{\alpha}(e^{i\phi(\cdot)}\,f)(x)= e^{i\phi(x)}\,D_{\alpha} f(x) +e^{i\phi(x)} \,\Lambda_{\alpha}((e^{i(\phi(x+y)-\phi(x))}-1)f(\cdot+y))(x).
 \end{equation}

 Now we assume that
 \begin{equation}
\label{d5}
n=1,\hskip10pt \alpha\in(0,1),\hskip10pt \phi(x)=\phi_{t}(x)=t x^3,
\end{equation}
we shall obtain a bound for
\begin{equation}
\label {goal1}
\begin{aligned}
&\|\Lambda_{\alpha}((e^{i(\phi(x+y)-\phi(x))}-1)f(\cdot+y))\|_p\\
\\
&=\|\,\lim_{\epsilon \to 0} \int_{|y| \geq \epsilon} \frac{e^{i(\phi(x+y)-\phi(x))}-1}{|y|^{1+\alpha}}\,f(x+y) dy\,\|_p.
\end{aligned}
\end{equation}

 We restrict ourselves to the case $\alpha\in (0,1)$ which allows us to perform estimates by passing the absolute value inside the integral sign in
 \eqref{d3}.

  We recall the elementary estimates
  \begin{equation}
  \begin{aligned}
  \label{d6}
  \begin{cases}
  &(a)\;\;\;\;\forall \,\theta\in \R\;\;\;\;\;\;|e^{i\theta}-1|\leq 2,\\
  &(b) \;\;\;\;\forall \,\theta\in\R\;\;\;\;\;\;|e^{i\theta}-1|\leq 2|\sin(\theta/2)|\leq |\theta|.
  \end{cases}
  \end{aligned}
  \end{equation}

  Combining \eqref{d6} (a) and Minkowski's integral inequality  it follows that
   \begin{equation}
  \label{step1}
  \begin{aligned}
 &  \|\, \int_{|y| \geq 1/100} \frac{e^{i(\phi(x+y)-\phi(x))}-1}{|y|^{1+\alpha}}\,f(x+y) \,dy\,\|_p \\
 \\
  & \leq \, \int_{|y| \geq 1/100} \frac{2}{\;\;|y|^{1+\alpha}}\,\| f(\cdot+y)\|_p \,dy\, \leq c_{\alpha} \|f\|_p.
  \end{aligned}
  \end{equation}

  So, it remains to estimate
   \begin{equation}
\label {goal2}
\|\,\lim_{\epsilon \to 0} \int_{\epsilon\leq |y| \leq 1/100 } \frac{e^{i(\phi(x+y)-\phi(x))}-1}{|y|^{1+\alpha}}\,f(x+y) dy\,\|_p.
\end{equation}

From  \eqref{d6} (b) and the mean value theorem one has that
\begin{equation}
  \label{step2}
  |e^{i(\phi(x+y)-\phi(x))}-1|\leq |\phi(x+y)-\phi(x)|= |y| \,|\int^1_0 \phi'(x+sy)ds|,
  \end{equation}
  with
  \begin{equation}
  \label{mvt}
  \phi'(x)=3t x^2.
  \end{equation}

  In particular, if $\,|x|\leq 100$ one has
  $$
    |e^{i(\phi(x+y)-\phi(x))}-1|\leq c\,|t|\,|y|,
    $$
   and
\begin{equation}
  \label{step3}
  \begin{aligned}
& \|\,\lim_{\epsilon \to 0} \int_{\epsilon\leq |y| \leq 1/100 } \frac{e^{i(\phi(x+y)-\phi(x))}-1}{|y|^{1+\alpha}}\,f(x+y) dy\,\|_{L^p(B_{100}(0))}
 \\
 \\
 &\leq \,|t|\,\int_{|y|\leq 1/100}\,
\frac{\| f(\cdot+y) \|_{L^p(B_{100}(0))}}
  {|y|^{\alpha}}\,dy
 \leq \,c_{\alpha} \,|t|\,\|f\|_p.
 \end{aligned}
 \end{equation}

  From the above estimates we can restrict ourselves in \eqref{goal1} to the case:
  $$
  |y|\leq 1/100,\;\;\;\;\;\text{and}\;\;\;\;\;|x|\geq 100.
  $$

  We sub-divide it into two parts:
 \begin{equation}
  \label{2sets}
  (a)\;\;\;|y|\,|x|^2\leq 1,\;\;\;\;\;(b)\;\;\;|y|\,|x|^2\geq 1.
  \end{equation}

  In the case (a) in \eqref{2sets} we change variable, $\tilde y=|x|^2y$, use \eqref{d6} part (b), \eqref{step2}, \eqref{mvt},  Minkowski's inequality
  and a second change of variable to obtain the bound

  \begin{equation}
  \label{step4}
  \begin{aligned}
&  \|\,\int_{|y|\leq 1/|x|^{2}}\,\frac{ |t| |x|^{2}\,|f(x+y)|}{|y|^{\alpha}}dy\|_{L^p(\{|x|\geq100\})}\\
\\
&= \|\,\int_{|\tilde y|\leq 1}\,\frac{|t| |x|^{2\alpha}\,|f(x+\frac{\tilde y}{|x|^{2}})|}{|\tilde y|^{\alpha}}d\tilde y\|_{L^p(\{|x|\geq100\})}\\
\\
&\leq \|\,\int_{|\tilde y|\leq 1}\,\frac{|t| |x + \frac{\tilde y}{|x|^{2}}|^{2\alpha}\,|f(x+\frac{\tilde y}{|x|^2})|}{|\tilde y|^{\alpha}}d\tilde y\|_{L^p(\{|x|\geq100\})}
\\
\\
&\;\;\;+
\|\,\int_{|\tilde y|\leq 1}\,\frac{|t||\frac{\tilde y}{|x|^{2}}|^{2\alpha}\,|f(x+\frac{\tilde y}{|x|^{2}})|}{|\tilde y|^{\alpha}}d\tilde y
\|_{L^p(\{|x|\geq100\})}\\
\\
&\leq \,c_{\alpha} |t| (\| \,|x|^{2\alpha}\,f\|_p+\|f\|_p),
\end{aligned}
\end{equation}
since
\begin{equation}
\label{changeofvariable}
\frac{\tilde y}{|x|^{2}}=y,\;\;\;\;\; |y|\leq 1/100,\;\;\;\;\;|x|\geq 100,\;\;\;\; \;\text{so} \;\;\;d(x+\frac{\tilde y}{|x|^{2}})\sim dx.
 \end{equation}

  In the case (b) in \eqref{2sets} changing variable, $\tilde y=x^2y$, using \eqref{d6} part (a), Minkowski's inequality, and a second change of
  variable
  as in \eqref{changeofvariable} we get
   \begin{equation}
  \label{step5}
  \begin{aligned}
  &\|\,\int_{1/x^2\leq|y|\leq 1/100} \,\frac{|f(x+y)|}{|y|^{1+\alpha}}\,dy\|_{L^p(\{|x|\geq100\})}\\
  \\
  &=
  \| \,\int_{1\leq|\tilde y|\leq x^2/100}\,\frac{\;|x|^{2\alpha}}{|\tilde y|^{1+\alpha}} \,|f(x+\frac{\tilde y}{|x|^2})|\,d \tilde y\|_{L^p(\{|x|\geq100\})}
 \\
 \\
 &\leq c_{\alpha}  \int_{1\leq |\tilde y|} \| \,|x|^{2\alpha} \,f(x+\frac{\tilde y}{|x|^{2}})\,\chi_{\{|x|\geq 10|\tilde y|^{1/2}\}}(x)\| _{L^p(\{|x|\geq100\})} \,\frac{d \tilde y}{|\tilde
 y|^{1+\alpha}}\\
  \\
  &\leq c_{\alpha}  \int_{1\leq |\tilde y|} \| \, |x+\frac{\tilde y}{|x|^{2}}|^{2\alpha} f(x+\frac{\tilde y}{|x|^{2}})\chi_{\{|x|\geq 10|\tilde y|^{1/2}\}}(x)\|
  _{L^p(\{|x|\geq100\})}\,\frac{d \tilde y}{|\tilde y|^{1+\alpha}}\\
 \\
 &\;\;\; +   c_{\alpha}  \int_{1\leq |\tilde y|} \,\| f(x+\frac{\tilde y}{|x|^{2}})\,\,\chi_{\{|x|\geq 10|\tilde y|^{1/2}\}}(x)\| _{L^p(\{|x|\geq100\})}\,
 \frac{d \tilde y}{|\tilde y|^{1+\alpha}}\\
\\
 & \leq c_{\alpha}\, (\|f\|_p+ \|\,|x|^{2\alpha}\,f\|_p).
  \end{aligned}
  \end{equation}

  Therefore, collecting the above results we have  the proof of the following:
   \vskip.1in
 \begin{lemma}
 \label{Lemma1}
Let $n=1$, $\,\alpha\in (0,1)$,  and $p\in(1,\infty)$. If
$$
f\in  L^{\alpha,p}(\R)\cap L^p(|x|^{2\alpha p}dx),
$$
then for all $t\in \R$ and for almost every $x\in \R$

\begin{equation}
\label{d10}
D_{\alpha} (e^{itx^3}\,f)(x) = e^{itx^3}\, D_{\alpha}f(x) + e^{itx^3} \,\Phi_{t,\alpha}(f)(x),
\end{equation}
 with
\begin{equation}
\label{d10a}
\Phi_{t,\alpha}(f)(x)=\lim_{\epsilon \to 0} \frac{1}{c_{\alpha}}\int_{|y|\geq \epsilon} \frac{e^{it((x+y)^3
-x^3)}-1}{|y|^{1+\alpha}}\,f(x+y)dy,
\end{equation}

\begin{equation}
\label{d10-norm}
\| \Phi_{t,\alpha}(f)\|_p\leq c_{\alpha}(1+|t|)( \|f\|_p + \| \,|x|^{2\alpha} \,f\|_p),
\end{equation}
and $c_{\alpha}$ as in \eqref{pro1}.
\end{lemma}

 From the  proof of Lemma \ref{Lemma1} it follows that under appropriate assumptions on the regularity and the growth
 of  a symbol $\varphi :\R^n\to \R$ one has  that
$$
D_{j,\alpha} (e^{it\varphi(x)}\,f)(x) = e^{it\varphi(x)}\, D_{j,\alpha}f(x) + e^{it\varphi(x)} \Phi_{j,\varphi,t,\alpha}(f)(x),
$$
 with
\begin{equation}
 \label{def11}
 D_{j,\alpha}f(x) = \lim_{\epsilon \to 0} \frac{1}{c_{\alpha}} \,\int_{|y_j| \geq \epsilon}  \frac{f(x+y_j\,\vec e_j)-f(x)}{|y_j|^{1+\alpha}}dy_j,
 \end{equation}
$$
\Phi_{j,\varphi,t,\alpha}(f)(x)=\lim_{\epsilon \to 0} \frac{1}{c_{\alpha}}\int_{|y_j|\geq \epsilon}
\frac{e^{it(\varphi(x+y_j\,\vec e_j)-\varphi(x))}-1}{|y_j|^{1+\alpha}}\,f(x+y_j\,\vec e_j)dy_j,
$$
and
 $$
\| \Phi_{j,\varphi,t,\alpha}(f)\|_p\leq c_{\alpha}(1+|t|)( \|f\|_p + \| |\partial_{x_j}\varphi(x)|^{\alpha} \,f\|_p),
$$
for $j=1,\dots,n$.

  Next, we consider the unitary group of operators $\{U(t) : t\in\R\}$ in $L^2(\R)$ defined as
  \begin{equation}
  \label{def33}
  U(t)u_0(x)=U(t)u_0(x)= ( e^{it \xi^3}\widehat{u}_0(\xi))^{\lor}(x).
  \end{equation}

  Thus, for $\alpha\in(0,1)$ using \eqref{pro1} one has that

  $$
  |x|^{\alpha}\,U(t)u_0(x)= |x|^{\alpha}( e^{it \xi^3}\widehat{u}_0(\xi))^{\lor}(x)= (D_{\alpha} (e^{it \xi^3}\widehat{u}_0(\xi)))^{\lor}(x).
  $$
 and from Lemma \ref{Lemma1} that \begin{equation}
  \label{for1}
  D_{\alpha} (e^{it\xi^3}\,\widehat u_0)(\xi) = e^{it\xi^3}\, D_{\alpha}\widehat u_0(\xi) + e^{it\xi^3} \Phi_{t,\alpha}(\widehat
  u_0)(\xi),
  \end{equation}
   with
  $$
  \| \Phi_{t,\alpha}(\widehat u_0)\|_p\leq c_{\alpha}(1+|t|)( \|\widehat u_0\|_p + \| \,|\xi|^{2\alpha} \,\widehat u_0\|_p).
  $$

  Hence, taking Fourier transform in \eqref{for1}
  we obtain the identity
 \begin{equation}
  \label{for2}
 |x|^{\alpha} \,U(t) u_0(x)= U(t)(|x|^{\alpha} u_0)(x) + U(t) (\{ \Phi_{t,\alpha}(\widehat u_0)(\xi)\}^{\lor})(x).
  \end{equation}
   with $\Phi_{t,\alpha}$ as in \eqref{d10a} and
  \begin{equation}
  \label{estimate1}
 \begin{aligned}
 & \| \{ \Phi_{t,\alpha}(\widehat u_0)(\xi)\}^{\lor}\|_2=\|\Phi_{t,\alpha}(\widehat u_0)\|_2
 \\
 & \leq c_{\alpha}(1+|t|) ( \|\widehat u_0\|_2 + \| \,|\xi|^{2\alpha} \,\widehat u_0\|_2)\\
 & \leq c_{\alpha}(1+|t|) (\|u_0\|_2 + \|\,D^{2\alpha}u_0\|_2).
 \end{aligned}
  \end{equation}
    Moreover, we claim that if $\beta\in(0,\alpha)$, then
\begin{equation}
  \label{for3}
D^{\beta}( |x|^{\alpha} \,U(t) u_0)(x)= U(t)(D^{\beta}|x|^{\alpha} u_0)(x) + U(t) (D^{\beta}\{ \Phi_{t,\alpha}(\widehat u_0)(\xi)\}^{\lor})(x).
  \end{equation}
  with
 \begin{equation}
  \label{estimate2}
  \|\,D^{\beta} (\{ \Phi_{t,\alpha}(\widehat u_0)(\xi)\}^{\lor})(x)\|_2\leq c_{\alpha,\beta}(1+|t|)(\|u_0\|_2 + \|\,D^{\beta+2\alpha}u_0\|_2).
  \end{equation}

 Notice that for $u_0\in \mathcal S(\R)$ the identities  \eqref{for2} and \eqref{for3} hold pointwise for each $(x,t)\in \R^2$.

  To prove \eqref{estimate2} we need to show that
  \begin{equation}
  \label{est7}
  \begin{aligned}
 &  \| D^{\beta}_x
(\int \frac{e^{it((\xi+\eta)^3-\xi^3)}-1}
 {|\eta|^{1+\alpha}}\,\widehat u_0(\xi+\eta)\, d\eta)^{\lor} \|_2
\\
\\
 &\leq c_{\alpha,\beta}(1+|t|) (\|u_0\|_2+\|D^{\beta+2\alpha}u_0\|_2).
  \end{aligned}
  \end{equation}
 Thus, we write
   \begin{equation}
  \label{est8}
  \begin{aligned}
 &  \| D^{\beta}_{x}(\int\frac{e^{it((\xi+\eta)^3-\xi^3)}-1}{|\eta|^{1+\alpha}}\,\widehat u_0(\xi+\eta)\,
 d\eta)^{\lor} \|_2\\
& = \| \int\frac{ |\xi|^{\beta}(e^{it((\xi+\eta)^3-\xi^3)}-1)}{|\eta|^{1+\alpha}}\,\widehat u_0(\xi+\eta)\,
 d\eta \|_2\\
& \leq  \| \int\frac{ |\xi|^{\beta}|e^{it((\xi+\eta)^3-\xi^3)}-1|}{|\eta|^{1+\alpha}}\,|\widehat u_0(\xi+\eta)|\,
 d\eta \|_2\\
 & \leq c_{\beta} \| \int\frac{ |\xi+\eta|^{\beta}|e^{it((\xi+\eta)^3-\xi^3)}-1|}{|\eta|^{1+\alpha}}\,|\widehat u_0(\xi+\eta)|\,
 d\eta \|_2\\
&\; \;\;+c_{\beta}  \| \int\frac{ |\eta|^{\beta}|e^{it((\xi+\eta)^3-\xi^3)}-1|}{|\eta|^{1+\alpha}}\,|\widehat u_0(\xi+\eta)|\,
 d\eta \|_2\\
 &= \Omega_1+\Omega_2.
  \end{aligned}
 \end{equation}

   Following the argument used in the proof of Lemma \ref{Lemma1} to get \eqref{d10-norm} one has that
   \begin{equation}
  \label{est9}
  \begin{aligned}
 &  \Omega_1\leq c_{\alpha}(1+|t|)(\| |\xi|^{\beta} \widehat u_0\|_2
 +\|\,|\xi|^{2\alpha}|\xi|^{\beta} \widehat u_0\|_2)\\
  & \;\;\;\;=c_{\alpha}(1+|t|)(\| D^{\beta}u_0\|_2+\| D^{\beta+2\alpha}u_0\|_2).
  \end{aligned}
 \end{equation}

 To bound $\Omega_2$  we observe that its estimate is similar to that used in the proof of Lemma \ref{Lemma1}
 with $\alpha-\beta$ instead of $\alpha$. Hence,
 \begin{equation}
  \label{est10}
  \begin{aligned}
 &  \Omega_2\leq c_{\alpha}(1+|t|)(\|  \widehat u_0\|_2+\|\,|\xi|^{2(\alpha-\beta)} \widehat u_0\|_2)\\
  & \;\;\;\;=c_{\alpha}(1+|t|)(\| u_0\|_2+\| D^{2(\alpha-\beta)}u_0\|_2).
  \end{aligned}
 \end{equation}

   Collecting the above information one obtains the proof of Theorem \ref{point-wise}.

 \underline{Remarks:} a) From the proof of Theorem 1, it is clear that  \eqref{for2}--\eqref{estimate2}
 hold for $f(\cdot,t)$ instead of $u_0$ with the suitable modifications.
 
  (b) The hypothesis $\beta \in(0,\alpha)$ in Theorem \ref{point-wise} is necessary to bound

 $$
 \aligned
 &\| D^{\beta}(\{ \Phi_{t,\alpha}(\widehat u_0)(\xi)\}^{\lor})\|_2=\| |\xi|^{\beta}\Phi_{t,\alpha}(\widehat u_0)(\xi)\|_2,
 \\
 \\
 =& \| \int\frac{ |\xi|^{\beta}(e^{it((\xi+\eta)^3-\xi^3)}-1)}{|\eta|^{1+\alpha}}\,\widehat u_0(\xi+\eta)\,
 d\eta \|_2\\
 \endaligned
 $$
  in the region where  $ |\xi+\eta|\leq |\xi|/10$ with $|\xi|\sim |\eta|\gg 1$.

 (c) We observe that if $u_0\in\mathcal S(\R)$, then the pointwise  identities  \eqref{d20}-\eqref{d21}  hold for all $(x,t)\in \R^2$.
  Therefore a density argument and  the Strichartz estimate associated to the group $\{U(t)\}$ (see \cite{KePoVe91})
  \begin{equation}
  \label{str}
 ( \int_{-\infty}^{\infty}\| U(t)u_0\|_{\infty}^{6}dt)^{1/6}\leq c\|u_0\|_2,
 \end{equation}
 show that under the hypotheses of Theorem \ref{point-wise},  \eqref{d20}-\eqref{d21} hold for all $x\in\R$ almost everywhere
  $t\in\R$.
  
  (c) The result in Theorem \ref{point-wise} also holds and the proof is similar to the above one for solutions of the linear IVP,
  \begin{equation}\label{DGBO}
\begin{cases}
\partial_t u -   D^{1+a}_x\partial_x u = 0, \qquad t, x\in \R,\;\;\;\; 0\leq a < 1,\\
u(x,0) = u_0(x),
\end{cases}
\end{equation}
where $D^s$ denotes  the homogeneous derivative of order $s\in\R$,
$$
D^s=(-\partial_x^2)^{s/2}\;\;\;\text{so}\;\;\;D^s f=c_s\big(|\xi|^s\widehat{f}\,\big)^{\vee}, \;\;\;\text{with} \;\;\;D^s=(\mathcal H\,\partial_x)^s,
$$
and $\mathcal H$ denotes the Hilbert transform,
\begin{equation*}
\mathcal H f(x)=\frac{1}{\pi}\lim_{\epsilon\downarrow 0}\int\limits_{|y|\ge \epsilon} \frac{f(x-y)}{y}\,dy=(-i\,\sgn(\xi) \widehat{f}(\xi))^{\vee}(x).
\end{equation*}

It is not clear how to employ Theorem \ref{point-wise}  to obtain solutions via contraction  of the IVP associated to the equation above with nonlinear
term like the one of the KdV equation called dispersion generalized Benjamin-Ono (DGBO) equation. Nevertheless there are optimal persistency results in weighted 
Sobolev spaces via energy estimates for the solutions of the IVP associated to the DGBO equation \cite{FLP}.

\section{Proofs of Theorem \ref{th4} and Theorem \ref{th5}}\label{S:2}

\underline{Proof of Theorem \ref{th4}}: 

We shall restrict our attention to the most interesting case $s=1/4$ and $r=1/8$, i.e. $u_0\in Z_{1/4,1/8}$.

We begin with  a brief review of  the argument used in the proof of Theorem A  in \cite{KPV1}. The details of this proof will be used later to complete  the proof of Theorem \ref{th4}. 

 First, let us assume that
 $$
u_0\in \dot H^{1/4}(\R).
$$

 For $w:\R\times[-T,T]\to\R$ with $T $ to be fixed below, define
\begin{equation}
\label{norm1}
\begin{aligned}
\mu_1^T(w) =& \|D_x^{1/4}w\|_{L^{\infty}_TL^2_x}+ \|\partial_x w\|_{L_x^{20} L_T^{5/2}}  
+ \|D_x^{1/4} w\|_{L_x^{5} L_T^{10}}\\
&+\|D_x^{1/4}\; \partial_x w\|_{L_x^\infty L_T^2}+ \|w\|_{L_x^{4} L_T^\infty}.
\end{aligned}
\end{equation}

 Denote by $\Phi(v)=\Phi_{u_0}(v)$
the solution of the linear inhomogeneous IVP
\begin{equation}
\label{inhoeq}
\partial_tu+\partial^3_x u +v^2{\partial_x v}=0, \hskip10pt
u(x,0)=u_0(x).
\end{equation}
The idea is to apply the contraction principle  to  the integral equation version of the IVP \eqref{inhoeq}, i.e.
\begin{equation}
\label{in-inhoeq}
u(t)=\Phi(v(t))=U(t)u_0-\int_0^{t} U(t-t')(v^2 \; \partial_x v)
(t')dt'.
\end{equation}

From  the linear estimates concerning the group $\{U(t)\,:\,t\in\R\}$ established in \cite{KPV1} one has that
\begin{equation}
\label{est1a}
\mu_1^T(U(t)u_0)\leq c_0 \| D^{1/4}_xu_0\|_2,\;\;\;\;\;\;\;\;\forall\,T>0.
\end{equation}

Here and below $c_0$ will denote a universal constant whose value  may change (increase) from line to line.
 Hence,
\begin{equation}
\label{est1b}
\begin{aligned}
&\mu_1^T(\int_0^t U(t-t') v^2\partial_x v(t')dt')\\
\\
& \leq c_0 \| D_x^{1/4}(v^2\partial_x v)\|_{L^1_TL^2_x}\leq
c_0 T^{1/2} \| D_x^{1/4}(v^2\partial_x v)\|_{L^2_xL^2_T}.
\end{aligned}
\end{equation}

  Using the  calculus of inequalities in the Appendix in \cite{KPV1} (Theorem A.8)
one gets  that
\begin{equation}
\label{ineq1}
\begin{aligned}
&\|D_x^{1/4}(v^2 \; {\partial_x v})\|_{L_x^2 L_T^2} \\
&\leq c_0\|D_x^{1/4}(v^2)\|_{L_x^{20/9} L_T^{10}}\|\partial_xv
\|_{L_x^{20} L_T^{5/2}}\;+\;
c_0\|v^2\|_{L_x^2 L_T^\infty}\|D_x^{1/4} \; {\partial_x v}
\|_{L_x^\infty L_T^2} \\
&\leq c_0\|v\|_{L_x^4 L_T^\infty}\|D_x^{1/4} v\|_{L_x^5 L_T^{10}} \; \|\partial_xv
\|_{L_x^{20} L_T^{5/2}}
\; + \;  c_0\|v\|^2_{L_x^4 L_T^\infty}\|D_x^{1/4} \; {\partial_x v}\|_{L_x^\infty L_T^2}\\
&\leq c_0(\mu_1^T(v))^3.
\end{aligned}
\end{equation}

Inserting the estimates \eqref{est1a}, \eqref{est1b}, and \eqref{ineq1} in the integral equation \eqref{in-inhoeq} it follows
that
\begin{equation}
\label{ineq3}
\begin{aligned}
\mu_1^T(\Phi(v))&\leq c_0\|D^{1/4}_xu_0\|_2+c_0\,\int_0^T \|D_x^{1/4}(v^2 \; {\partial_x v})
\|_2(t)dt \\
&\leq c_0 \|D^{1/4}_xu_0\|_2 +c_0\,T^{1/2}(\mu_1^T(v))^3.
\end{aligned}
\end{equation}

A similar argument leads to the estimate 
\begin{equation}
\label{ineq4}
\mu_1^T(\Phi(v)-\Phi(\widetilde v))\leq c_0\,T^{1/2}(\mu_1^T(v)+\mu_1^T(\widetilde v))^2\,
\mu_1^T(v-\widetilde v).
\end{equation}

This basically proves the main part of Theorem A. More precisely, one has that the operator $\Phi=\Phi_{u_0}$ in \eqref{in-inhoeq} defines a contraction in the set
\begin{equation}
\label{set1}
\{v\ :\R\times [-T, T]\to \R\,:\,\mu_1^T(v)\leq 2c_0\|D^{1/4}_xu_0\|_2\},
\end{equation}
with 
\begin{equation}
\label{time-est}
T=\frac{1}{32 \,c_0^6\,\|D^{1/4}_xu_0\|_2^4}.
\end{equation}

Hence, the IVP \eqref{kgKdV} with $k=2$ has a unique solution $u\in C([-T,T] :\dot H^{1/4}(\R))$
satisfying 
\begin{equation}
\label{est-no}
\mu_1^T(u)\leq 2c_0\|D^{1/4}_xu_0\|_2,
\end{equation}
with $T$ as in \eqref{time-est}.

 Now, we assume that
 $$
 u_0\in H^{1/4}(\R),
 $$
 and define
$$
\mu_2^{T_0}(w) = \|w\|_{L^{\infty}_{T_0}L^2_x}+  \|\partial_x w\|_{L_x^{\infty} L_{T_0}^{2}} + \|w\|_{L^{6}_{T_0}L^{\infty}_x} + \mu_1^{T_0}(w),
$$
with $\mu_1^{T_0}$ defined in \eqref{norm1} and $\,T_0>0$ to be fixed below. By the previous argument we have a solution $u=u(t)$ in the class defined by \eqref{norm1} of the integral equation 
\begin{equation}
\label{eq}
u(t)=U(t)u_0-\int_0^{t} U(t-t')(u^2 \; \partial_x u)
(t')dt'.
\end{equation}
By \eqref{est1a} and Strichartz estimates \eqref{str} one has that
\begin{equation}
\label{aaa1}
\| U(t)u_0\|_{L^{\infty}_{T_0}L^2_x} +\|\,\partial_x U(t)u_0\|_{L^{\infty}_xL^2_{T_0}} + \| U(t)u_0\|_{L^{6}_{T_0}L^{\infty}_x} \leq c_0\|u_0\|_2,\;\;\forall\;T_0>0.
\end{equation}
 
 Therefore
\begin{equation}
\label{aaa2}
\begin{aligned}
&\|\,\int_0^t U(t-t')u^2\partial_xu(t')dt'\,\|_{L^{\infty}_{T_0}L^2_x} 
+ \|\,\partial_x \int_0^t U(t-t')u^2\partial_xu(t')dt'\,\|_{L^{\infty}_xL^2_{T_0}}\\
&+ \|\,\int_0^t U(t-t')u^2\partial_xu(t')dt'\,\|_{L^{6}_{T_0}L^{\infty}_x} \\
&\;\;\leq c_0\,\|u^2\partial_x u\|_{L^1_{T_0}L^2_x}\leq c_0\,T_0^{1/2}\,\|u^2\partial_x u\|_{L^2_xL^2_{T_0}}\\
&\;\;\leq c_0\,T_0^{1/2}\,\|u^2\|_{L^2_xL^{\infty}_{T_0}}\,\|\partial_x u\|_{L^{\infty}_xL^2_{T_0}}\leq c_0
\,T_0^{1/2}\,\|u\|^2_{L^4_xL^{\infty}_{T_0}}\,\|\partial_x u\|_{L^{\infty}_xL^2_{T_0}}\\
&\;\;\leq c_0
\,T_0^{1/2} (\mu_1^{T_0}(u))^2 \mu_2^{T_0}(u).
\end{aligned}
\end{equation}

Collecting the above information one has that
$$
 \mu_2^{T_0}(u) \leq 2c_0(\|u_0\|_2+\|D^{1/4}_xu_0\|_2) + c_0
\,T_0^{1/2} (\mu_1^{T_0}(u))^2 \mu_2^{T_0}(u).
$$
Hence, taking 
 $T_0=T$ as in \eqref{time-est}, i.e. 
\begin{equation}
\label{key2}
c_0T^{1/2} (\mu_1^{T}(u))^2\leq 1/2,
\end{equation}
 it follows
\begin{equation}
\label{007}
\mu_2^{T}(u)\leq 4 c_0(\|u_0\|_2+\|D^{1/4}_xu_0\|_2).
\end{equation}

By uniqueness we have 
 $$
 u\in C([-T,T]\,:\, H^{1/4}(\R))\cap L^6([-T,T]:L^{\infty}(\R)) \cap\dots
$$
which can be extended to the interval $ [-T^*,T^*] $ as far as the 
\begin{equation}
\label{key1}
\sup_{t\in[-T^*,T^*]} \|D^{1/4}u(t)\|_2<\infty,
\end{equation}
since we recall that the $L^2$-norm of the real solutions of the IVP \eqref{kgKdV} is preserved in time. 
Now we turn our attention to the most interesting case in Theorem \ref{th4}
 $$
 u_0\in Z_{1/4,1/8}=H^{1/4}(\R)\cap L^2(|x|^{1/4}\,dx),
 $$
 and introduce the notation
 $$
\mu_3^{T_0}(w) =\mu_2^{ T_0}(w) + \|\,|x|^{1/8}w(t)\|_{L^{\infty}_{ T_0}L^2_x},
$$
with $\,T_0>0$ to be fixed below.

From Theorem \ref{point-wise}  (see \eqref{d20}-\eqref{d20-norm}) and the linear estimates in \eqref{aaa1} it follows that
\begin{equation}
\label{est1c}
\mu_3^{T_0}(U(t)u_0)\leq c_0 \| |x|^{1/8} u_0 \|_2+ c_0 (1+{T_0})
(\|u_0\|_2+ \| D_x^{1/4}u_0\|_2).
\end{equation} 

Now taking $\varphi \in C^{\infty}_0(\mathbb R)$ with $\varphi =1,\;|x|<1/2$ and $\,\varphi=0,\;|x|\geq1$ we write
\begin{equation}\label{neww}
\begin{aligned}
|x|^{1/8}u^2\partial_xu &= \varphi(x) |x|^{1/8}u^2\partial_xu + (1-\varphi(x))|x|^{1/8}u^2\partial_xu \\
=&\; \varphi |x|^{1/8}u^2\partial_xu + \partial_x((1-\varphi)|x|^{1/8}u^3/3)-\partial_x((1-\varphi)|x|^{1/8})u^3/3\\
\equiv&\;  A_1+A_2+A_3.
\end{aligned}
\end{equation}

Same argument as in \eqref{est1c}  and \eqref{neww} yield
\begin{equation}\label{nonlinear}
\begin{split}
&\||x|^{1/8}\,\int_0^t U(t-t') u\partial_x u(t')dt'\|_{L^2_x} \\
&\le \|\int_0^t U(t-t')(A_1+A_2+A_3)\,dt'\|_{L^2_x}
+\int_0^T \| \{\Phi_{t,1/4}\widehat{(u\partial_x u)} \}^{\vee}\|_{L^2_x}\,dt'\\
&\le \int_0^{T_0} \|U(t-t')(A_1+A_3)\|_{L^2_x}dt' +\| \int_0^t U(t-t') A_2\,dt'\|_{L^2_x}\\
&\;\;\;+ c_0(1+T_0)\int_0^{T_0} (\|u^2\partial_x u\|_{L^2_x}
+\|D^{1/4}(u^2\partial_x u)\|_{L^2_x})\,dt.
\end{split}
\end{equation}

Thus,
\begin{equation}\label{est1d}
\begin{aligned}
 \int_0^{T_0} \|U(t-t')A_1\|_{L^2_x}dt' &\leq c_0\,\|u^2\partial_x u\|_{L^1_{T_0}L^2_x}\leq c_0\,T_0^{1/2}\,\|u^2\partial_x u\|_{L^2_xL^2_{T_0}}\\
&\leq c_0\,T_0^{1/2}\,\|u^2\|_{L^2_xL^{\infty}_{T_0}}\,\|\partial_x u\|_{L^{\infty}_xL^2_{T_0}}\\
&\leq c_0
\,T_0^{1/2}\,\|u\|^2_{L^4_xL^{\infty}_{T_0}}\,\|\partial_x u\|_{L^{\infty}_xL^2_{T_0}}.
\end{aligned}
\end{equation}
Using a duality argument (see \cite{KPV1}) one has that
$$
\|\partial_x \int_0^t U(t-t') F(t')dt'\|_{L^{\infty}_TL^2_x}\leq c \| F\|_{L^{1}_xL^2_T}.
$$
Hence,
\begin{equation}\label{est1dd}
\begin{aligned}
&\| \int_0^t U(t-t') A_2\,dt'\|_{L^2_x}\leq c_0 \| \,(1-\varphi)|x|^{1/8} u^3\|_{L^1_xL_{T_0}^2}\\
& \leq c_0\| \,|x|^{1/8} u\|_{L^2_xL^2_{T_0}}\,\|u^2\|_{L^2_xL_{T_0}^{\infty}}\leq c_0 T_0^{1/2} \| \,|x|^{1/8} u\|_{L^{\infty}_{T_0}L^2_x}\,\|u\|^2_{L^4_xL_{T_0}^{\infty}}.
\end{aligned}
\end{equation}

Finally,
\begin{equation}\label{est1ddd}
\begin{aligned}
& \int_0^{T_0} \|U(t-t')A_3\|_{L^2_x}dt' \leq c_0 \| \,u^3\|_{L^1_{T_0}L_x^2}\\
& \leq c_0\| \,u\|_{L^{\infty}_{T_0}L^4_x}^2 \|u\|_{L^1_{T_0}L^{\infty}_x}\leq c_0 T_0^{5/6} \| D^{1/4}_x u\|_{L^{\infty}_{T_0}L^2_x}^2 \|u\|_{L^6_{T_0}L^{\infty}_x}.
\end{aligned}
\end{equation}

Inserting the estimates \eqref{est1c}-\eqref{est1ddd}, \eqref{ineq1}, \eqref{est1b} and \eqref{aaa2} in the integral equation 
\eqref{eq} it follows
that
\begin{equation}
\label{ineq4b}
\begin{aligned}
\mu_3^{T_0}(u)& \leq c_0\| |x|^{1/8}u_0\|_2 + c_0(1+T_0)
(\|u_0\|_2+\|D^{1/4}_xu_0\|_2)\\
&\;\;\;+c_0 T_0^{1/2}(\mu_1^{T_0}(u))^2 \mu_3^{T_0}(u)\\
&\;\;\; + c_0(1+T_0^{1/3})T_0^{1/2}(\mu_1^{T_0}(u))^2\,\mu_2^{T_0}(u).
\end{aligned}
\end{equation}

Thus, taking $T_0\in(0,T]$ with $T$ as in \eqref{key2} and \eqref{007} one can rewrite \eqref{ineq4b}  as
\begin{equation}
\label{ineq4bbb}
\begin{aligned}
\mu_3^{T_0}(u)& \leq  2\,c_0\| |x|^{1/8}u_0\|_2 + 2 c_0(1+T_0)
(\|u_0\|_2+\|D^{1/4}_xu_0\|_2)\\
&\;\;\; + 4 c_0 (1+T_0^{1/3})\mu_2^{T_0}(u)\\
&\leq 2\,c_0\| |x|^{1/8}u_0\|_2 + 20 c_0(1+T_0) (\|u_0\|_2+\|D^{1/4}_xu_0\|_2)
\end{aligned}
\end{equation}
which basically completes the proof of Theorem \ref{th4}.

\vskip3mm

\underline{Proof of Theorem \ref{th5}}: 

 We shall consider the most interesting case $s=1/4$, and recall that the $L^2$-norm of the solution $u(t)$ is preserved.  
 
 By Theorem B for any given $T^*>0$ and $u_0\in H^{1/4}(\R)$ one has that the corresponding solution 
 $u=u(x,t)$ of the IVP \eqref{kgKdV} with $k=2$ satisfies 
 $$
 u\in C([-T^*,T^*] \,:\,H^{1/4}(\R))\cap....
 $$
with
 $$
 K=\max_{[-T^*,T^*]}\|D_x^{1/4}u(t)\|_2.
 $$
 Following \eqref{time-est} we define
 $$
 T'=\frac{1}{32 \,c_0^{6}\,K^4},
 $$
and split the interval $[-T^*,T^*]$ into $\;2T^*/T'$ sub-intervals. In each of these sub-intervals we can apply Theorem \ref{th4}  observing that the right hand side of \ref{ineq4bbb} depends on $K$, $\;2T^*/T'$ and the initial value $\,\| |x|^{1/8}u_0\|_2$  to get the desired solution to the whole interval $[-T^*,T^*]$.  
 
\vskip3mm

{\underbar{ACKNOWLEDGMENTS}}: F.L was partially supported by CNPq and FAPERJ-\/Brazil. G.P. was supported by NSF grant DMS-1101499.

\end{document}